\documentclass[10pt]{article}
\usepackage{graphicx}
\usepackage{amsmath}
\usepackage[dvips]{epsfig}
\usepackage{amssymb}

\makeatletter
\renewcommand{\section}{\@startsection
  {section}%
  {2}%
  {0mm}%
  {\baselineskip}%
  {0.3 \baselineskip}%
  {\centering}}
\makeatother
\begin{document}

\title { On the alternating sums of powers of
consecutive $q$-integers }
\author{ Taekyun Kim\\[0.5cm]
  Institute of Science Education,\\
        Kongju National University, Kongju 314-701,  Korea\\
          {\it e-mail: tkim $@$kongju.ac.kr }  }

\date{}
\maketitle

{\footnotesize {\bf Abstract}\hspace{1mm}
  In this paper we construct $q$-Genocchi numbers and polynomials.
  By using these numbers and polynomials, we investigate the
  $q$-analogue of alternating sums of powers of consecutive
  integers due to Euler.



\section{Introduction }
The Genocchi numbers $G_m$
 are defined by the generating function:
 $$F(t)=\dfrac{2t}{e^t+1}=\sum_{m=0}^\infty G_m \dfrac{t^m}{m!}, ( |t|< \pi ), \text{ cf. [3]} \eqno(1)$$
 where we use the technique method notation by replacing $G^m$ by
 $G_m ( m \geq 0)$ symbolically.
 It satisfies $G_1=1, G_3=G_5=G_7= \cdots=0, $ and even
 coefficients are given $G_m=2(1-2^{2m})B_{2m}=2mE_{2m-1},$
 where $B_m$ are Bernoulli numbers and $E_m$ are Euler numbers
 which are defined by
 $$\dfrac{2}{e^t+1}=\sum_{m=0}^\infty E_m \dfrac{t^m}{m!}, \text{ cf. [5, 6] }.$$
 For $x \in \mathbb{R}$ (= the field of real numbers), we consider
 the Genocchi polynomials as follows:
 $$F(x,t)=F(t)e^{xt}=\dfrac{2t}{e^t+1}e^{xt}=\sum_{n=0}^\infty G_n(x) \dfrac{t^n}{n!}. \eqno(2)$$
Note that $G_m(x)=\sum_{k=0}^m \binom mk G_k x^{m-k} $. Let us
also define the Genocchi polynomials of order $r$ as follows:
$$2 \left(\dfrac{1}{1+e^t}\right)^r e^{xt}=\sum_{n=0}^\infty G_n^{(r)}(x) \dfrac{t^n}{n!}, \text{ cf. [3]}.$$
In the special case $x=0$, we define $G_n^{(r)}(0)=G_n^{(r)}$.
What is the value of the following sum for a given positive
integer $k$ ?
$$1^k+2^k+3^k+ \cdots +n^k.$$
Let us denote this sum by $f_k(n)$. Finding formulas for $f_k(n)$
has interested mathematicians for more than 300 years since the
time of Jacob Bernoulli, cf. [1, 7, 9].  It was well known that
$$f_n{(k-1)}=\frac{1}{n+1} \sum_{i=0}^n \binom{n+1}iB_i k^{n+1-i}, \text{ cf. [9] },\eqno(3) $$
where $\binom nk $ is binomial coefficients.

Let $n, k$ be positive integers ($k>1$), and let
$$T_n(k)=-1^k+2^k-3^k+4^k-5^k+\cdots+(-1)^{k-1}(n-1)^k.$$
Following an idea due to Euler, it was known that
$$T_n(k)=\frac{(-1)^{k+1}}{2}\sum_{l=0}^{k-1}\binom nl E_lk^{n-l}+
\frac{E_n}{2}\left(1+(-1)^{k+1}\right),\mbox{ cf. [5].} \eqno(4)
$$

Let $q$ be an indeterminate which can be considered in complex
number field, and for any integer $k$ define the $q$-integer as
$$[k]_q=\dfrac{q^k-1}{q-1}=1+q+ \cdots + q^{k-1}.$$
Throughout this paper we assume that  $q \in \mathbb{C} $ with
$0<q<1$.  Recently many authors studied $q$-analogue of the sums
of powers of consecutive integers.  In [2], Garrett and Hummel
gave a combinatorial proof of a $q$-analogue of $\sum_{k=1}^n k^3=
\binom{n+1}k ^2 $ as follows:
$$\sum_{k=1}^n q^{k-1}[k]_q^2 \left( \left[ \dfrac{k-1}{2}\right]_{q^2}+
 \left[\dfrac{k+1}{2}\right]_{q^2}\right)
=\left[ \begin{array}{c} n+1 \\  2 \end{array} \right]_q^2, $$
where
$$\left[ \begin{array}{c} n \\  k \end{array} \right]_q=
\prod_{j=1}^k \dfrac{[n+1-j]_q}{[j]_q} \mbox{ denotes the }
q\mbox{-binomial coefficients}.
$$
Garrett and Hummel, in their paper, asked for a simpler
$q$-analogue of the sums of cubes.  As a response to Garrett and
Hummel 's question, Warnaar gave a simple $q$-analogue of the sums
of cubes as follows:
$$\sum_{k=1}^n q^{2n-2k}[k]_q^2[k]_{q^2}
=\left[ \begin{array}{c} n+1 \\  2 \end{array} \right]_q^2 ,
\mbox{ cf. [10] }.
$$
Let $$f_{m,q}(n)=\sum_{k=1}^n [k]_{q^2}[k]_q^{m-1}
q^{(n-k)\frac{m+1}{2}}.$$ Then we note that $\lim_{q \to 1}
f_{m,q}(n)=f_m(n) $.

Warnaar [10] ( for $m=3$) and Schlosser [8] gave formulae for
$m=1,2,3,4,5 $ as the meaning of the $q$-analogues of the sums of
consecutive integers, squares, cubes, quarts and quints.  Let $n,
k ( >1) $ be the positive integers. In the recent paper, it was
known that
$$ \sum_{j=0}^{k-1} q^j [j]_q^n=\dfrac{1}{n+1} \sum_{j=0}^n \binom{n+1}j
\beta_{j,q} q^{kj} [k]_q^{n+1-j} - \dfrac{1}{n+1}(1-q^{(n+1)k})
\beta_{n+1,q}, \mbox{ see [4] }, $$ where $\beta_{j,q}$ are called
Carlitz's  $q$-Bernoulli numbers. Originally $q$-Genocchi numbers
and polynomials were introduced by Kim-Jang-Pak in 2001 [3], but
they do not seem to be the most natural ones. In this paper we
give another construction of a $q$-Genocchi numbers and
polynomials which are different than $q$-Genocchi numbers and
polynomials of Kim-Jang-Pak in 2001 [3].
  By using these numbers and polynomials, we investigate the
  $q$-analogue of alternating sums of powers of consecutive
  integers.

\section{ $q$-Genocchi numbers and polynomials }
Let $F_{q,k}(t)$ be the generating functions of the $q$-Genocchi
numbers as follows:
$$F_{q,k}(t)=[2]_q t \sum_{j=0}^\infty
q^{k-j}[j]_{q^2} (-1)^{j-1}\exp \left( t[j]_{q^2}
q^{\frac{k-j}{2}}\right) =\sum_{j=0}^\infty G_{n,q}
\dfrac{t^n}{n!}
$$
By using Taylor expansion in the above, we see that
$$\aligned
\sum_{j=0}^\infty G_{n,k,q} \dfrac{t^n}{n!} &=[2]_q t
\sum_{j=0}^\infty q^{k-j}[j]_{q^2} (-1)^{j-1} \sum_{n=0}^\infty
\dfrac{[j]_q^n q^{\frac{n(k-j)}{2}}}{n!} t^n \\
&=[2]_q t \sum_{j=0}^\infty q^{k-j}[j]_{q^2} (-1)^{j-1}
\sum_{n=0}^\infty \left \{ \dfrac{1}{(1-q)^n} q^{\frac
{n(k-j)}{2}}
\sum_{m=0}^n \binom nm (-1)^m q^{jm} \right \} \dfrac{t^n}{n!} \\
&=\dfrac{[2]_q t}{1-q^2} \sum_{n=0}^\infty \dfrac{1}{(1-q)^n}
q^{k+ \frac{nk}{2}} \sum_{m=0}^n \binom nm (-1)^m  \left(
\sum_{j=0}^\infty (-1)^{j-1} q^{mj-j -\frac{nj}{2}} (1-q^{2j})
\right)  \dfrac{t^n}{n!} \\
&=\dfrac{[2]_q t}{1-q^2 } \sum_{n=0}^\infty \dfrac{1}{(1-q)^n}
\sum_{m=0}^{n} \binom nm (-1)^m
\dfrac{q^{m-1-\frac{n}{2}+k+\frac{nk}{2}}(1-q^2)}
{(1+q^{m-1-\frac{n}{2}})(1+q^{m+1-\frac{n}{2}})} \dfrac{t^n}{n!}.
\endaligned
$$
Note that $G_{0,k,q}=0 $. Hence, we have
$$\aligned
\sum_{n=1}^\infty G_{n,k,q} \dfrac{t^n}{n!} & ={ t}
\sum_{n=1}^\infty \dfrac{1}{(1-q)^{n}} \sum_{m=1}^{n} \binom
{n-1}{m-1} \dfrac{(-1)^{m-1} q^{m+k
+\frac{(n-1)k}{2}-\frac{n-1}{2}-2}}
{(1+q^{-2+m-\frac{n-1}{2}})(1+q^{m-\frac{n-1}{2}})}
\dfrac{t^{n-1}}{(n-1)!}\\
&=  \sum_{n=1}^\infty \dfrac{1}{(1-q)^{n}}\sum_{m=1}^{n} \binom
{n}{m} \dfrac{(-1)^{m-1}m q^{m+k +\frac{(n-1)(k-1)}{2}-2}}
{(1+q^{-2+m-\frac{n-1}{2}})(1+q^{m-\frac{n-1}{2}})}
\dfrac{t^{n}}{n!} \endaligned
$$
By comparing the coefficients of $ \dfrac{t^n}{n!} $ on  both
sides of the above equation, we obtain the below:
\bigskip

 {\bf Theorem 1.}\quad
 Let $k, n( n \geq 1)$ be positive integers. Then we have
 $$ G_{n,k,q}=  \left( \dfrac{1}{1-q} \right )^{n}\sum_{m=1}^n \binom nm
  \dfrac{(-1)^{m-1}m q^{m+k +\frac{(n-1)(k-1)}{2}-2}}
{(1+q^{-2+m-\frac{n-1}{2}})(1+q^{m-\frac{n-1}{2}})} .$$

We also define the generating function.
$$F_{q,k}(t,k)=[2]_q t \sum_{j=0}^\infty
q^{-j}[j+k]_{q^2} (-1)^{j+k-1}\exp \left( t[j+k]_{q}
q^{-\frac{j}{2}}\right) =\sum_{n=0}^\infty G_{n,k,q}(k)
\dfrac{t^n}{n!}
$$
By using the binomial theorem and some elementary calculations in
the above equation, we have
$$\aligned
\sum_{n=0}^\infty G_{n,k,q}(k) \dfrac{t^n}{n!} &=\dfrac{[2]_q
t}{1-q^2}  \sum_{j=0}^\infty q^{-j}(1-q^{2j+2k}) (-1)^{j+k-1}
\sum_{n=0}^\infty [j+k]_q^n q^{-\frac{j}{2}n} \dfrac{t^n}{n!}\\
&=\dfrac{[2]_q t}{1-q^2}  \sum_{j=0}^\infty q^{-j}(1-q^{2j+2k})
(-1)^{j+k}
\sum_{n=0}^\infty \left\{ \left( \dfrac{1}{1-q}\right)^{n} \sum_{m=0}^n \binom nm
(-1)^m q^{jm+km}  q^{-\frac{j}{2}n} \right \} \dfrac{t^n}{n!}\\
&=\dfrac{[2]_q t}{1-q^2} \sum_{n=0}^\infty  \left(
\dfrac{1}{1-q}\right)^{n} \sum_{m=0}^n \binom nm (-1)^{m+k} q^{mk}
 \sum_{j=0}^\infty (1-q^{2j+2k})
(-1)^{j}q^{-j+jm-\frac{n}{2}j} \dfrac{t^n}{n!} \\
&=\dfrac{[2]_q t}{1-q^2}\sum_{n=0}^\infty \left(
\dfrac{1}{1-q}\right)^{n} \sum_{m=0}^n \binom nm (-1)^{m+k} q^{mk}
\left \{ \dfrac{1}{1+q^{m-\frac{n}{2}-1}}
 -\dfrac{q^{2k}}{1+q^{1+m-\frac{n}{2}}}  \right \} \dfrac{t^n}{n!} \\
&=t \sum_{n=0}^\infty  \left( \dfrac{1}{1-q}\right)^{n+1}
\sum_{m=0}^n \binom nm (-1)^{m+k} \left(
\dfrac{q^{mk}}{1+q^{m-\frac{n}{2}-1}}
 -\dfrac{q^{(m+2)k}}{1+q^{1+m-\frac{n}{2}}} \right)\dfrac{t^n}{n!}  \\
&=t \sum_{n=1}^\infty  \left( \dfrac{1}{1-q}\right)^{n}
\sum_{m=1}^n \binom {n-1}{m-1} (-1)^{m-1+k} \left(
\dfrac{q^{(m-1)k}}{1+q^{m- 2-\frac{n-1}{2}}}
 -\dfrac{q^{(m+1)k}}{1+q^{m-\frac{n-1}{2}}}
 \right)\dfrac{t^{n-1}}{(n-1)!}\\
 &= \sum_{n=1}^\infty  \left( \dfrac{1}{1-q}\right)^{n}
\sum_{m=1}^n \binom {n}{m} (-1)^{m-1+k} \left( \dfrac{m
q^{(m-1)k}}{1+q^{m- 2-\frac{n-1}{2}}}
 -\dfrac{m q^{(m+1)k}}{1+q^{m-\frac{n-1}{2}}}
 \right)\dfrac{t^{n}}{n!}.
 \endaligned
$$
Note that $G_{0, k,q}(k)=0$.  Therefore we obtain the following
theorem.

\bigskip

 {\bf Theorem 2.}\quad
 Let $k,  n( n \geq 1)$ be positive integers. Then we have
 $$ G_{n,k,q}(k)=\left( \dfrac{1}{1-q}\right)^{n}
\sum_{m=1}^n \binom {n}{m} (-1)^{m-1+k} \left( \dfrac{m
q^{(m-1)k}}{1+q^{m- 2-\frac{n-1}{2}}}
 -\dfrac{m q^{(m+1)k}}{1+q^{m-\frac{n-1}{2}}}
 \right).$$

\bigskip

 {\bf Remark 3.} Note that $$
 (1)  \lim_{q\to 1} G_{n,k,q}=G_n^{(2)}, \text{ }(2)  \lim_{q\to 1} G_{n,k,q}(k)\neq G_n^{(2)}(k).
$$

It is easy to see that
$$\aligned
&[2]_q t  \sum_{j=0}^\infty q^{k-j}[j]_{q^2}(-1)^{j-1} \exp \left(
t[j]_{q^2} q^{\frac{k-j}{2}} \right)  - [2]_q t \sum_{j=0}^\infty
q^{-j}[j+k]_{q^2}(-1)^{j-1+k} \exp \left(
t[j+k]_{q} q^{-\frac{j}{2}} \right) \\
&=[2]_q t  \sum_{j=0}^{k-1}(-1)^{j-1}[j]_{q^2}q^{k-j}  \exp \left(
t[j]_{q} q^{\frac{k-j}{2}} \right).
 \endaligned
$$
Thus, we easily see that
$$[2]_q   \sum_{j=0}^{k-1}[j]_{q^2} (-1)^{j-1}[j]_{q}^{n-1}
q^{\frac{(k-j)(n+1)}{2}} = \dfrac{G_{n,k,q}-G_{n,k,q}(k)}{n}.$$
Therefore we obtain the following theorem.

\bigskip

{\bf Theorem 4.}\quad
 Let $k,  n( n \geq 1)$ be positive integers. Then we have
 $$   \sum_{j=0}^{k-1}[j]_{q^2} (-1)^{j-1}[j]_{q}^{n-1}
q^{\frac{(k-j)(n+1)}{2}} =
\dfrac{G_{n,k,q}-G_{n,k,q}(k)}{n[2]_q}$$

\end{document}